\documentclass{ifacconf}

\usepackage{graphicx}      
\usepackage{natbib}        
\usepackage{amsmath} 
\usepackage{amsfonts}
\usepackage{amssymb}
\usepackage{enumerate}
\usepackage{graphicx}
\usepackage{epstopdf}
\usepackage{subfigure}
\usepackage{url} 
\usepackage{tikz}
\usepackage{color}
\usepackage{float}

\newtheorem{lemma}{Lemma}
\newtheorem{theorem}{Theorem}
\newtheorem{corollary}{Corollary}
\begin{document}
\begin{frontmatter}

\title{Population Control of {\it Giardia lamblia}
} 


\author[First]{Victor~Hugo~Pereira~Rodrigues} 
\author[Second]{Maria~Fantinatti}
\author[First]{Tiago~Roux~Oliveira}
\author[First]{Wilton~dos~Santos~Freitas}

\address[First]{\mbox{Department of Electronics and Telecommunication Engineering (DETEL),} \\ State University of Rio de Janeiro (UERJ), Rio de Janeiro, Brazil. \\ \mbox{\!\!\!\!\!\!\!\!\!\!\!\!\!\!\!\!  (e-mails: rodrigues.vhp@gmail.com, tiagoroux@uerj.br, wiltonmessias55@gmail.com)}}
\address[Second]{\mbox{Interdisciplinary Medical Research Laboratory,} \\ Oswaldo Cruz Foundation, Rio de Janeiro, RJ, Brazil. \\ (e-mail: fantinatti@ioc.fiocruz.br)}

\begin{abstract}                
Giardia lamblia is a flagellate intestinal protozoan with global distribution causing the disease known as giardiasis. This parasite is responsable for 35.1\% of outbreaks of diarrhea caused by contaminated water which and mainly affects children in whom it can cause physical and cognitive impairment. In this paper, we consider a model of population dynamics 
to represent the behavior of Giardia lamblia in vitro, taking into account its mutation characteristic that guarantees to the protozoan resistance to the drug metronidazole. Different from what is found in the literature, it is pursued as the control objective the extermination of the protozoan considering that the parameters of the model are uncertain and only the partial measurement of the state vector is possible. On these assumptions, a control law is designed and the stability of the closed-loop system is rigorously proved. Simulation and experimental results illustrate the benefits of the proposed population control method of {\it Giardia lamblia}.
\end{abstract}


\end{frontmatter}

\section{Introduction}
Giardia is a single-celled organism responsible for giardiasis. This protozoan is divided into species according to morphological and molecular characteristics, and the host it infects. The species {\it Giardia duodenalis} (syn. {\it G. lamblia, G. intestinalis}) is responsible for the disease in mammals, including domestic animals and humans.

The parasite {\it G. duodenalis} has two evolutionary forms: cyst and trophozoite. The infection starts by ingestion of cysts present in water or on the contaminated food surfaces. This cyst, when passing through the stomach region, begins a process of excystation and in the first portions of the small intestine takes the trophozoite form. The trophozoites proliferate (by binary division) and colonize the duodenum. When carried by the intestinal flow, they initiate encystation and are released into the feces in the cyst form already infected \citep{CL:2010}.

Infection with {\it G. duodenalis} has a global distribution. The highest frequencies are observed in developing countries, such as Brazil (over 50\%), since their transmission is closely related to poor basic sanitation conditions \citep{GS:2002,JKC:2013}.

The most frequent symptoms of giardiasis are a consequence of acute or chronic diarrhea, such as abdominal colic, flatulence, dehydration, nausea, vomiting and fatigue \citep{TRM:2013,FX:2011}. Notice that the diarrhea and pneumonia are the mean causes of children's death worldwide \citep{DK:2004} and it is estimated that 35.1\% of diarrhea outbreaks by means of contaminated water are provided by {\it Giardia} infections \citep{BK:2011}. Moreover, prolonged exposure to the parasite can cause delayed growth and development, poor cognitive function and detrimental effects on nutritional status in children \citep{F:1996,BLGLB:2002}.

Metronidazole is the most used drug in the world for the giardiasis treatment and is the drug of first choice recommended by the Ministry of Health of Brazil. In most cases, treatment is effective, however, the number of cases of persistent infection after using the drug has been reported \citep{NLAC:2015,FOFTVBC:2020}. In this context, more studies need to be carried out in order to elucidate the real effects of the drug's action on the parasite, helping to understand the cases of refractory to treatment.

In \citep{LLZ:2013}, by using evolutionary game theory, the authors provides a model for biological systems subjects to dynamics that exhibit mutation behavior. Moreover, an optimal control strategy was developed and validated for drug delivery using the pathogen {\it G. lamblia}. Nevertheless, the design is based on the assumptions of the full-state measurement and assuming the perfect knowledge of the model parameters.

In this paper, a new control strategy is designed to eliminate/attenuate the flagellate {\it G. lamblia} population employing the mathematical model inspired by \citep{LLZ:2013} and assuming that only the plant output is available for feedback in a scenario where all model parameters are uncertain. After introducing the control law, the stability analysis of the proposed closed-loop system is carried out. Numerical simulations and experimental tests are also performed in order to evaluate the proposed control scheme for drug delivery.


\section{Problem Formulation}

The first dynamic model that describes the population behavior of {\it G. lamblia}, {\it in vitro}, was presented by \citep{LLZ:2013}. This approach takes in count the resistance of the delivery drug (metronidazol) due to parasite mutation. 

In this paper, the following growth model (\ref{eq:dotX1})--(\ref{eq:y}) is directly related to the time-varying growth rate $r(t)$ in (\ref{eq:r}) determined by the physiology of the {\it G. lamblia} and the carrying capacity rate $K$ determined by the {\it G. lamblia} characteristic and/or environmental factor: 
\begin{align}
		\dot{x}_{1}\left(t\right)&=r\left(t\right)x_1\left(t\right)\left(1-\frac{x_1\left(t\right)}{K}\right) \,, \label{eq:dotX1}\\
		\dot{x}_{2}\left(t\right)&=-\sigma^{2}\beta_{m}x_{2}\left(t\right)+\frac{\sigma^{2}\beta_{m}\sqrt{w_{m}}}{\sqrt{2}}\mbox{exp}\left(-\frac{x_{2}^{2}\left(t\right)}{w_{m}}\right)\frac{x_{1}\left(t\right)}{K} \,, \label{eq:dotX2} \\
	y\left(t\right)&=x_{1}\left(t\right)\,. \label{eq:y}
\end{align}
As shown in \citep{LLZ:2013}, the {\it G. lamblia} population $x_{1}(t)$ [Cells/ml] is driven by the time-varying growth rate  
\begin{align}
r\left(t\right)=r_{0}\frac{x_{1}(t)}{K}+\beta_{m}x_{2}\left(t\right)\mbox{exp}\left(-\frac{x_{2}^{2}\left(t\right)}{w_{m}}\right)-\beta_{d}u\left(t\right)\,, \label{eq:r}
\end{align}
where $r_{0}>0$ is the natural growth rate, representing the case of no mutation and no drug. The state variable $x_{2}(t) \in \mathbb{R}$ is the drug resistance by mutation while the $\beta_{m}>0$ is the resistance rate. The input signal is the metronidazol dosage $u(t) \in \mathbb{R}$ [ug/ml] and the constant $\beta_{d}>0$ represents the drug efficiency. The state is defined by $x (t) = [x_{1}(t)\,,x_{2}(t)]^{T} \in \mathbb{R}^{2}$ and the output is denoted by $y(t) \in \mathbb{R}$.  

Throughout the paper the following assumptions are considered:

\begin{description}
\item[{\bf (A1)}] All parameters of the plant (\ref{eq:dotX1})--(\ref{eq:r}) are uncertain. 
%
In practice, biological systems can vary significantly their parameters during the experimental process. Thus, the assumption (A1) is important to capture the degree of robustness that the closed-loop system must possess with respect to the parameter uncertainties. On the other hand, the designer should pay attention to avoid an overestimation that leads to higher control efforts. 
\item[{\bf (A2)}] Only the output $y(t)$ is available for the feedback control design. Due to financial and/or physical constrains, the full-state measurement is barely possible. For instance, the state variable $x_{2}(t)$ represents a mutation behavior of a protozoan population that obviously cannot be measured in practice.
\item[{\bf (A3)}] The control design must take in count the strictly positive nature of the problem, {\it i.e.}, $u(t)\geq 0$ for all time-instant $t$. Since the control signal is a delivered drug to an organism, one cannot consider negative control actions. 
\end{description}

In this sense, the control objective is the design of a robust output-feedback control law $u(t)$ for the nonlinear uncertain plant (\ref{eq:dotX1})--(\ref{eq:r})  that ensures, in closed-loop, the global stabilization. To overcome the lack of the full-state measurement, in next section a norm observer is introduced to upper bound $|x_{2}(t)|$  and allows for the design of the control law.

\section{Norm Observer for the {\it G. lamblia} Mutation}

The first-order approximation filter (FOAF) \citep{CCH:2008}, also called norm observer, is an important tool to upper bound the norm of an unmeasured state variable. The basic idea is to provide a norm bound for $|x_2(t)|$ using the available signals as well as the lower and upper bounds of the model parameters.

Lemma~\ref{lemma:normObserver} below shows how the norm observer gives us such an instantaneous norm bound for the unmeasured variable through the solution of 
\begin{align}
\dot{\hat{x}}_{2}(t)=-\lambda\hat{x}_{2}(t)+\gamma |y(t)|\,, \label{eq:hatX2_1}
\end{align}
with scalars $\lambda \in ~  (  0\,, \sigma^{2}\beta_{m} )$ and $\gamma \in ~ \left( \frac{\sigma^{2}\beta_{m}\sqrt{w_{m}}}{K\sqrt{2}} \,, +\infty \right) $. 
The state variable $\hat{x}_2(t)$ satisfies 
\begin{align}
|x_{2}(t)| \leq |\hat{x}_{2}(t)|+\pi(t)\,, \quad \forall t \geq 0\,, \label{eq:hatX2_2}
\end{align}  
where $\pi(t)$ is an exponentially decaying term depending on the initial conditions $x_{2}(0)$ and $\hat{x}_{2}(0)$,
with $x_{2}(t)$ provided in (\ref{eq:dotX2}).

\begin{lemma} \label{lemma:normObserver}
Consider the $x_{2}$-dynamics in (\ref{eq:dotX2}) and suppose that assumptions (A1)--(A3) are satisfied. Then, $\hat{x}_{2}(t)$ in (\ref{eq:hatX2_1}) is a norm observer of $x_{2}(t)$ satisfying (\ref{eq:hatX2_2}).
\end{lemma}

\begin{pf}
Consider the following candidate for the Lyapunov function,
\begin{align}
V_{1}(x_{2})=x_{2}^{2}(t)\,, \label{eq:V1_1}
\end{align}
whose time-derivative is 
\begin{align}
&\dot{V}_{1}(x_{2})\mathbb{=}2x_{2}(t)\dot{x}_{2}(t)\,, \label{eq:dotV1_1} \\
&\mathbb{=}2|x_{2}(t)|\frac{d|{x}_{2}(t)|}{dt}\,, \label{eq:dotV1_2} \\
&\mathbb{=}2x_{2}(t)\left[\mathbb{-}\sigma^{2}\beta_{m}x_{2}\left(t\right)\mathbb{+}\frac{\sigma^{2}\beta_{m}\sqrt{w_{m}}}{\sqrt{2}}\mbox{exp}\left(\!\!\mathbb{-}\frac{x_{2}^{2}\left(t\right)}{w_{m}}\!\right)\!\frac{x_{1}\left(t\right)}{K}\right], \label{eq:dotV1_3} \\
&\mathbb{=}2\left[\mathbb{-}\sigma^{2}\beta_{m}x_{2}^{2}\left(t\right)\mathbb{+}\frac{\sigma^{2}\beta_{m}\sqrt{w_{m}}}{\sqrt{2}}x_{2}\left(t\right)\mbox{exp}\left(\!\!\mathbb{-}\frac{x_{2}^{2}\left(t\right)}{w_{m}}\!\right)\!\frac{x_{1}\left(t\right)}{K}\right]. \label{eq:dotV1_4}
\end{align}
Since $\mbox{exp}\left(-\frac{x_{2}^{2}\left(t\right)}{w_{m}}\right)\in~ \rbrack 0\,, 1\rbrack$ and $\frac{x_{1}\left(t\right)}{K}\in \lbrack 0\,, 1\rbrack$, equation (\ref{eq:dotV1_4}) is upper bounded by
\begin{align}
&\dot{V}_{1}(x_{2})\leq2\left[-\sigma^{2}\beta_{m}x_{2}^{2}\left(t\right)+\frac{\sigma^{2}\beta_{m}\sqrt{w_{m}}}{K\sqrt{2}}|x_{1}\left(t\right)||x_{2}\left(t\right)|\right] \nonumber \\
&=2|x_{2}\left(t\right)|\left[-\sigma^{2}\beta_{m}|x_{2}\left(t\right)|+\frac{\sigma^{2}\beta_{m}\sqrt{w_{m}}}{K\sqrt{2}}|y\left(t\right)|\right]\,. \label{eq:dotV1_5}
\end{align}
From (\ref{eq:dotV1_2}) and (\ref{eq:dotV1_5}), with the scalars $\lambda \in ~  (  0\,, \sigma^{2}\beta_{m} )$ and $\gamma \in ~ \left( \frac{\sigma^{2}\beta_{m}\sqrt{w_{m}}}{K\sqrt{2}} \,, +\infty \right) $, it is possible to find the following upper bound: 
%
\begin{align}
\frac{d|{x}_{2}(t)|}{dt}&\leq-\sigma^{2}\beta_{m}|x_{2}\left(t\right)|+\frac{\sigma^{2}\beta_{m}\sqrt{w_{m}}}{K\sqrt{2}}|y\left(t\right)| \nonumber \\
&\leq-\lambda|x_{2}\left(t\right)|+\gamma|y\left(t\right)|\,. \label{eq:dotAbsX2}
\end{align}
Then, invoking the Comparison Lemma \citep[p. 102]{K:2002}, the solution $\hat{x}_{2}(t)$ of (\ref{eq:hatX2_1}) is an upper bound for $|x_{2}\left(t\right)|$ since, by subtracting (\ref{eq:hatX2_1}) from (\ref{eq:dotAbsX2}), one has
\begin{align}
\frac{d(|{x}_{2}(t)|-\hat{x}_{2}(t))}{dt}&\leq-\lambda(|x_{2}\left(t\right)|-\hat{x}_{2}(t))\,. \label{eq:dotAbsX3}
\end{align}
Therefore, one can write
\begin{align}
|{x}_{2}(t)|-\hat{x}_{2}(t)&\leq e^{-\lambda t}(|x_{2}\left(0\right)|-\hat{x}_{2}(0))\,, 
\end{align}
and, consequently, 
\begin{align}
|{x}_{2}(t)|&\leq \hat{x}_{2}(t)+e^{-\lambda t}(|x_{2}\left(0\right)|-\hat{x}_{2}(0)) \nonumber \\
&\leq |\hat{x}_{2}(t)|+\underbrace{e^{-\lambda t}(|x_{2}\left(0\right)|+|\hat{x}_{2}(0)|)}_{\pi(t)}\,, \label{eq:dotAbsX4}
\end{align}
which completes the proof. \hfill $\square$
\end{pf}

\section{Drug Dosage Control}
\label{sec:drugDosage}

Metronidazole is the most commonly drug used in the treatment of giardiasis worldwide \citep{GH:2001}. Clinical trials have employed conventional dosing (generally 250 mg/dose) two or three times daily (for five to ten days), in short course (for one to three days), or single-dose daily therapy (2.0 or 2.4 g/dose) \citep{ZMA:1997}. In the five- to ten-day regimens, one gets efficacy ranging from 60\% to 100\% in adult and pediatric patients, with a median efficacy in both groups of 92\%. For children given doses between 15 and 22.5 mg/kg/day for regimens of five to ten days, efficacy ranges from 80 to 100\%. On the other hand, single-dose and short-term treatments (a high dose given daily) have been proposed to improve adherence without sacrificing efficacy. They were used on adults and children. These regimens are generally less effective, particularly if only one dose of metronidazole is administered. The efficacy of single-dose therapy ranges from 36 to 60\% if the drug is given for one day, rises to 67 to 80\% if the drug is given for two days, and 93 to 100\% for three days of treatment.

Clinical resistance and treatment failure against {\it G. lamblia} infection to the antiparasitic drug metronidazole, as well as furazolidone, has been reported in up to 20\% of cases \citep{BPS:1988,F:1996} with recurrence rates as high as 90\% \citep{ZMA:1997}. Metronidazole resistance can be easily induced {\it in vitro} by trophozoites growing at a sublethal concentration of the drug \citep{BPS:1988}.

In this section, metronidazole drug delivery strategies will be employed to study the population dynamics of {\it G. lamblia} {\it in vitro}. 

The following control law is proposed in order to guarantee asymptotic convergence to zero of the output variable $y(t)$ (population of {\it G. lamblia}):
\begin{align}
u(t)=\frac{\bar{r}_{0}|y(t)|+\bar{\beta}_{m}|\hat{x}_{2}(t)|+\eta}{\underline{\beta}_{d}} \,, \label{eq:increaseDosage}
\end{align}
with scalars $\bar{r}_{0}>r_{0}/K$, $\bar{\beta}_{m}>\beta_{m}$, $\underline{\beta}_{d}<\beta_{d}$ and $\eta>\bar{\beta}_{m}(|x_{2}\left(0\right)|+|\hat{x}_{2}(0)|)$. 
This control law is able to obtain output exponential stabilization in such a way that
\begin{align}
y(t)&\leq y(0) e^{-\delta t} \,, \quad t\geq 0 \,, \label{eq:yUpperBound2}
\end{align}
where $y(0)$ and $\eta$ are positive constants. The main stability results are summarized in the next theorem.

\begin{theorem} \label{theorem:drugDelivery}
Consider the system (\ref{eq:dotX1})-(\ref{eq:r}), the control law in (\ref{eq:increaseDosage}), the norm-state estimator of {\it Giardia lamblia} (\ref{eq:hatX2_1})-(\ref{eq:hatX2_2}) and suppose assumptions (A1)--(A4) are satisfied. Then, the output signal (\ref{eq:y}) converges to zero exponentially such that (\ref{eq:yUpperBound2}) is verified and the state $x(t)$ is uniformly bounded. \label{theorem:increaseDosage}
\end{theorem}

\begin{pf}
Consider the candidate to Lyapunov function
\begin{align}
V_{2}(t)=y^{2}(t) \,, \label{eq:V2-1}
\end{align}
whose time derivative satisfies
\begin{align}
&\dot{V}_{2}(t)\mathbb{=}2y(t)\dot{y}(t) \label{eq:dotV2-0} \\
&\mathbb{\leq}2\left(\!\! \frac{r_{0} x_{1}(t)}{K}\mathbb{+}\beta_{m}|x_{2}(t)|e^{\mathbb{-}(x_{2}^{2}(t)/w_{m})}\mathbb{-}\beta_{d}u(t)\!\! \right)\! y^{2}(t)\! \left(\!\! 1\mathbb{-}\frac{y(t)}{K}\!\!\right)\!. \nonumber
\end{align}
From (A1), recalling that $e^{-x^2_2(t)/\omega_m}\leq 1$, and by defining the positive constants $\bar{r}_{0}>r_{0}/K$, $\bar{\beta}_{m}>\beta_{m}$ and $\underline{\beta}_{d}<\beta_{d}$,  it is possible to write
\begin{align}
\dot{V}_{2}(t)&\mathbb{\leq}2\left(\bar{r}_{0}|y(t)|\mathbb{+}\bar{\beta}_{m}|x_{2}(t)|\mathbb{-}\underline{\beta}_{ d}u(t)\right)y^{2}(t)\left(1\mathbb{-}\frac{y(t)}{K}\right) \,. \label{eq:dotV2-1}
\end{align}
Using the norm observer (\ref{eq:hatX2_2}), inequality (\ref{eq:dotV2-1}) is rewritten as
\begin{align}
\dot{V}_{2}(t)&\mathbb{\leq}2\left(\bar{r}_{0}|y(t)|\mathbb{+}\bar{\beta}_{m}\left(|\hat{x}_{2}(t)| \mathbb{+}\pi(t)\right)\mathbb{-}\underline{\beta}_{d}u(t)\right)\times \nonumber \\
&\quad \times y^{2}(t)\left(1\mathbb{-}\frac{y(t)}{K} \right) \,. \nonumber
\end{align}
Then, for any positive constant $\eta$ sufficiently large such that $\eta>\bar{\beta}_{m}|x_{2}\left(0\right)|+\bar{\beta}_{m}|\hat{x}_{2}(0)|$, {\it i.e.}, $\eta=\bar{\beta}_{m}|x_{2}\left(0\right)|+\bar{\beta}_{m}|\hat{x}_{2}(0)|+\delta$, where $\delta$ is a small positive constant, by using the control law (\ref{eq:increaseDosage}), one has
\begin{align}
\dot{V}_{2}(t)&\leq-2\delta y^{2}(t)\left(1-\frac{y(t)}{K}\right) \,. \label{eq:dotV2-2}
\end{align}
Comparing (\ref{eq:dotV2-2}) with (\ref{eq:dotV2-0}), it is clear that
\begin{align}
\dot{y}(t)&\leq-\delta y(t)\left(1-\frac{y(t)}{K}\right) \,. \label{eq:dotY2-0}
\end{align}
Therefore, applying the Comparison Lemma \citep[p. 102]{K:2002}, an upper bound $\bar{y}(t)$ for $y(t)$ is obtained by means ofthe dynamical system solution
\begin{align}
\dot{\bar{y}}(t)&=-\delta \bar{y(t)}\left(1-\frac{\bar{y(t)}}{K}\right) \,,\quad t\geq 0 \quad \mbox{with} \quad \bar{y}(0)=|y(0)| \,. \label{eq:dotBarY2-0}
\end{align}
The nonlinear differential equation (\ref{eq:dotBarY2-0}) is a Riccati equation whose solution is given by
\begin{align}
\bar{y}(t)&=\frac{y(0)}{1-\dfrac{y(0)}{K}[e^{-\delta t}-1]}e^{-\delta t}\,, \quad t\geq 0 \,. \label{eq:barY2-0}
\end{align}
Since $\bar{y}(t)$ is an upper bound for $y(t)$,
\begin{align}
y(t)&\leq\bar{y}(t) \nonumber \\
&=\quad\frac{y(0)}{1+\dfrac{y(0)}{K}[1-e^{-\delta t}]}e^{-\delta t} \nonumber \\
&\leq\quad y(0)e^{-\delta t} \,, \quad t\geq 0 \,, \nonumber
\end{align}
which completes the proof. \hfill $\square$
\end{pf}

The corollary below guarantees the exponential of the full-state vector $x(t)=[x_1(t)\,, x_2(t)]^T$. 

\begin{corollary}
Consider the plant (\ref{eq:dotX1})-(\ref{eq:r}), the control law in (\ref{eq:increaseDosage}), the norm-state estimator of {\it Giardia lamblia} (\ref{eq:hatX2_1})-(\ref{eq:hatX2_2}) and suppose assumptions (A1)--(A4) are satisfied. From Lemma~\ref{lemma:normObserver}, the upper bound $\hat{x}_{2}(t)$ for the unmeasured state variable $x_{2}(t)$, given by the solution of (\ref{eq:hatX2_1}), is input-to-state stable (ISS) with respect to $x_{1}(t)$. Consequently, from Theorem~\ref{theorem:drugDelivery}, in closed-loop, the state variablle $x_{1}(t)$ converges exponentially to zero. Then, $x_{2}(t)$ goes to zero exponentially as well.
\end{corollary}

\section{Numerical Simulations}

In this section, the simulation results obtained for the closed-loop system constituted by the plant (\ref{eq:dotX1})-(\ref{eq:r}) and the control law (\ref{eq:increaseDosage}). In all simulations, the values for the plant parameters are $r_0=0.179527$ per hour, $K=6.596 \times 10^{7}$ cells/ml , $\beta_{d}=0.00874109$, $\beta_{m}=0.04162605$, $w_{m}=45.8677$ and $\sigma=0.52947229225$. The initial conditions were chosen as: $x_{1}(0)=10^{5}$ cells/ml, $x_{2}(0)=4.80$ and $\hat{x}_{2}(0)=11.25$. To mitigate unnecessarily higher control efforts, the norm observer (\ref{eq:hatX2_1}) is implemented with parameters $\lambda=1.14 \times 10^{-2}$ and $\gamma=1.70 \times 10^{-9}$. The other control parameters are $\bar{r}_0=3 \times 10^{-09}$, $\bar{\beta}_m=0.05$, $\underline{\beta}_d=0.007$, $\delta=0.024$ and $\eta=\delta$.

\subsection{Population behavior in the absence of Metronidazole}

\begin{figure}[!ht]
	\centering
	\subfigure[{\it Giardia lamblia} Population. \label{fig:giardiaSimulation_ol_x1}]{\includegraphics[width=4.2cm]{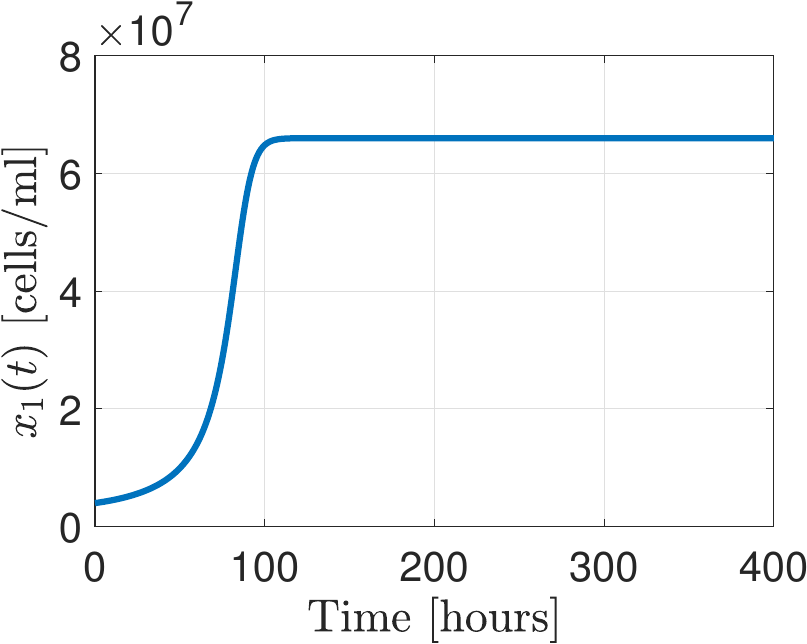}}
	~~
	\subfigure[{\it Giardia lamblia} Mutation. \label{fig:giardiaSimulation_ol_x2}]{\includegraphics[width=4.2cm]{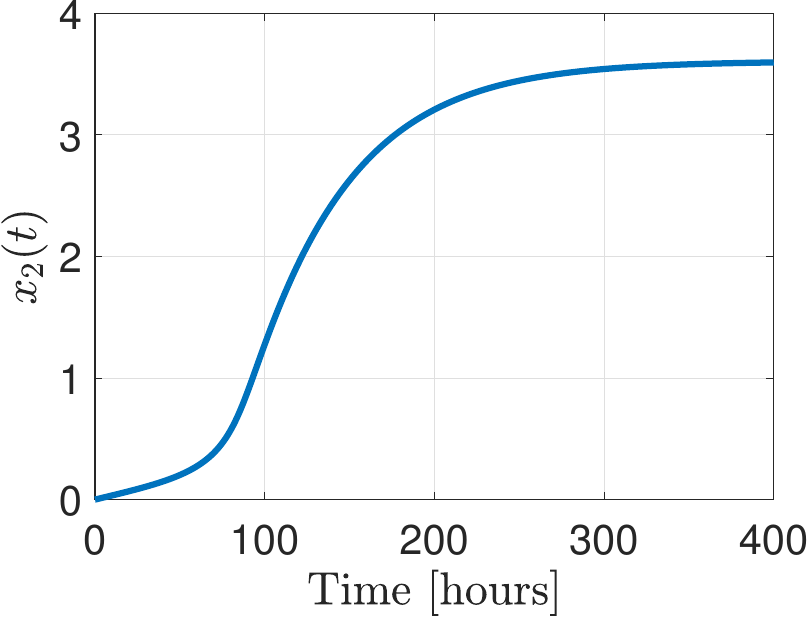}}
	\caption{Simulation results for the open-loop system.}
	\label{fig:simulation_ol}
\end{figure}

In this section, the simulation is conducted without any drug added to the culture. Figures~\ref{fig:giardiaSimulation_ol_x1} and \ref{fig:giardiaSimulation_ol_x2} show how both the trophozoites growth and their mutation variable increase monotonically without the drug delivered. After a certain time, the {\it Giardia Lamblia} population $y(t)=x_{1}(t)$ reaches the value of the carrying capacity $K$ and the state variable $x_{2}(t)$ tends to $\sqrt{\omega_{m}/2}$.  

\subsection{Population behavior in the presence of constant doses of Metronidazole}

Here, fixed-drug represents the drug dosage strategy where the drug dosage is the same at every time instant. It is easy to check that, if all the parameters of the plant (\ref{eq:dotX1})-(\ref{eq:r}) are known, it is possible to apply the control law
\begin{align}
u(t)=\frac{r_{0}+\beta_{m}\sqrt{\frac{\omega_{m}}{2}}+\delta}{\beta_{b}} \,, \quad  \forall t\geq 0\,, \label{eq:u2}
\end{align}
in closed-loop, ensuring that $y(t) \leq y(0)\exp(-\delta t)$, for all $t \geq 0$. 

\begin{figure}[!ht]
	\centering
	\subfigure[{\it Giardia lamblia} Population. \label{fig:giardiaSimulation_cte_yY0}]{\includegraphics[width=4.2cm]{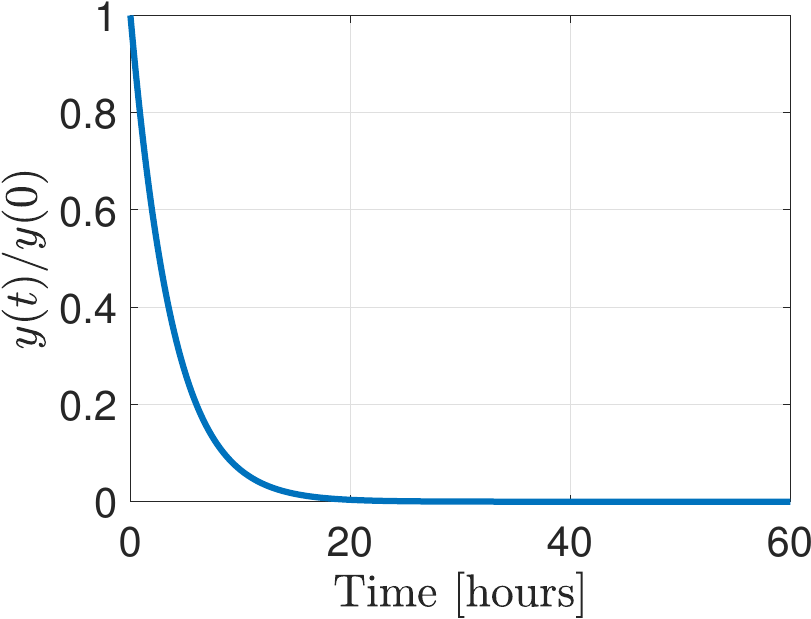}}
	~~
	\subfigure[Input signal; drug dosage in micro molar ($1\mu g/ml\approx 5,8425\mu M$). \label{fig:giardiaSimulation_cte_u}]{\includegraphics[width=4.2cm]{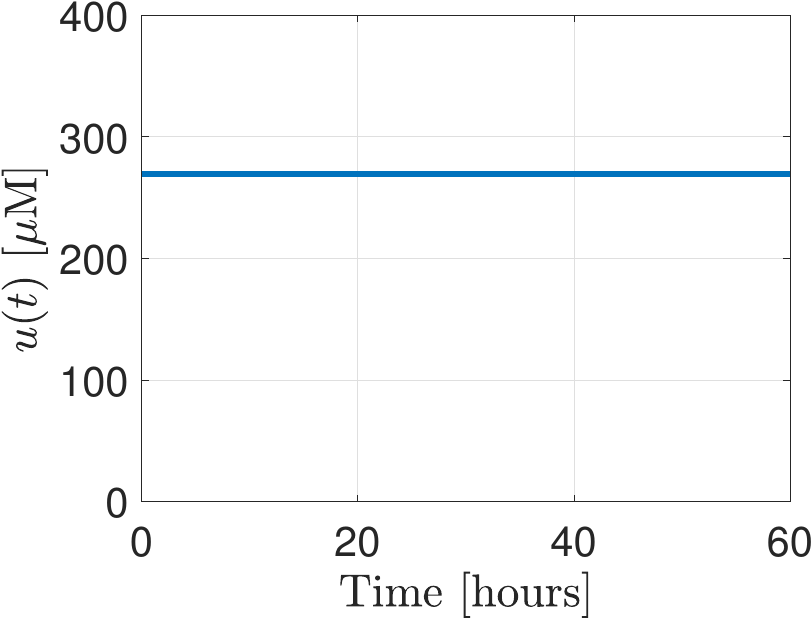}}
	\caption{Simulation results for the closed-loop system with a constant control law.}
	\begin{picture}(0,0)
		\put(-92,100.5){\includegraphics[width=2.75cm]{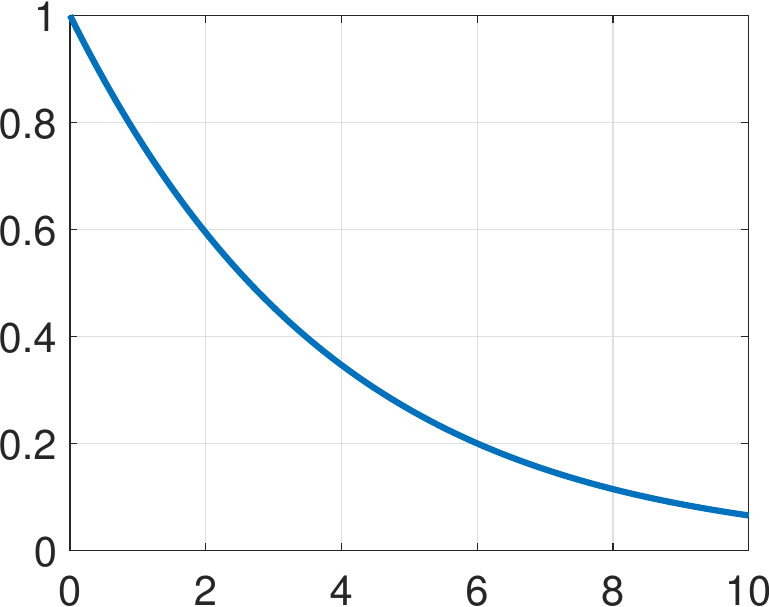}}
	\end{picture}
\end{figure}

In Figures~\ref{fig:giardiaSimulation_cte_yY0} and \ref{fig:giardiaSimulation_cte_u}, it can be seen that the population of {\it Giardia lamblia} decreases rapidly when the drug dosage is given by (\ref{eq:u2}).

\subsection{Population behavior in the presence of time-varying doses of Metronidazole}

\begin{figure}[!ht]
	\centering
	\subfigure[{\it Giardia lamblia} Population. \label{fig:giardiaSimulation_y}]{\includegraphics[width=4.2cm]{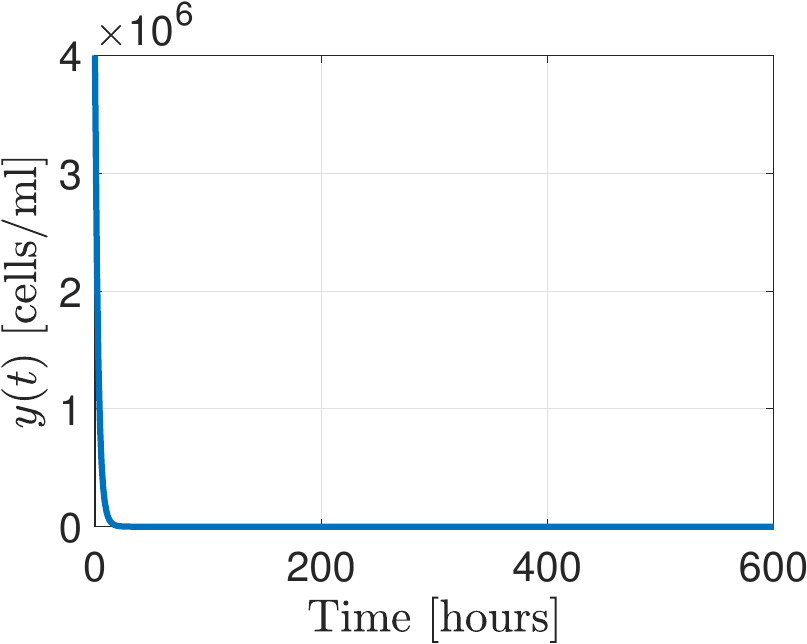}}
	~~
	\subfigure[{\it Giardia lamblia} Mutation and Norm Observer. \label{fig:giardiaSimulation_x2hatX2}]{\includegraphics[width=4.2cm]{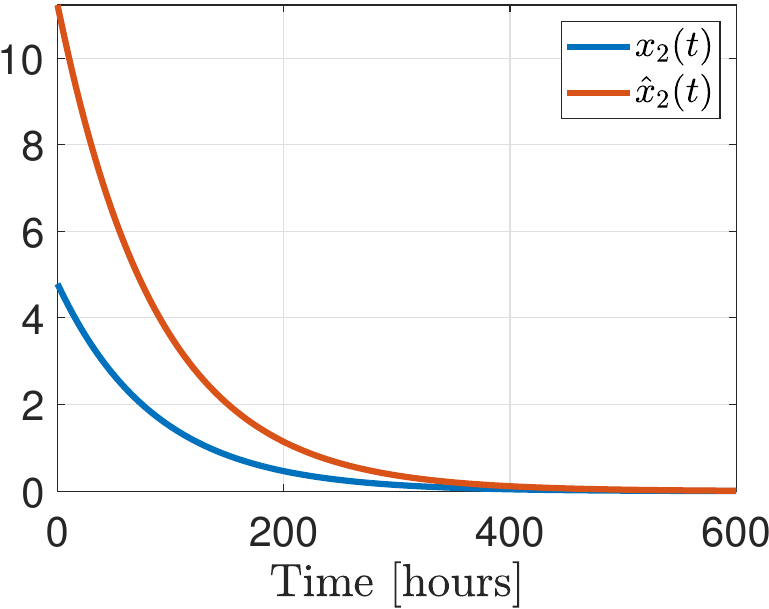}}
	\\
	\subfigure[Normalized population of {\it Giardia lamblia}. \label{fig:giardiaSimulation_yY0}]{\includegraphics[width=4.2cm]{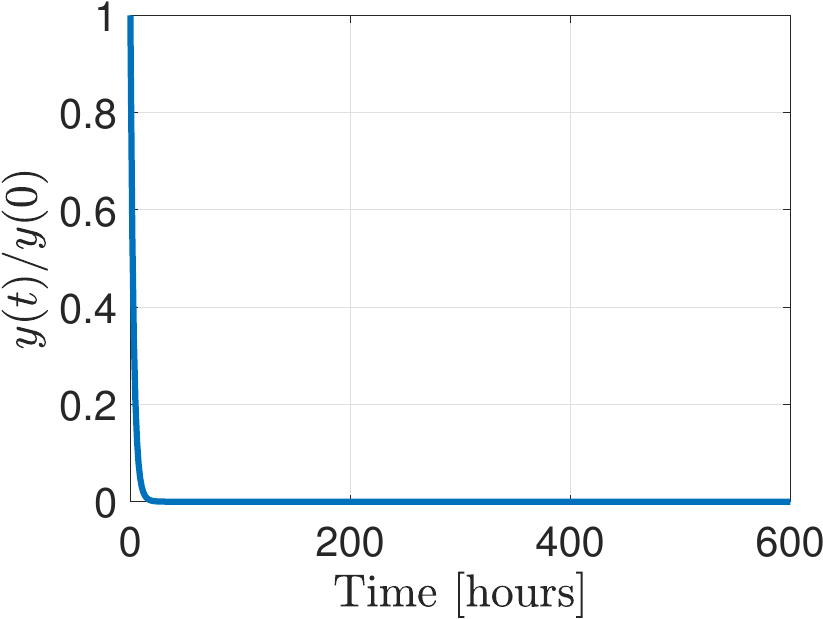}}
	~~
	\subfigure[Input signal; drug dosage in micro molar ($1\mu g/ml\approx 5,8425\mu M$). \label{fig:giardiaSimulation_u}]{\includegraphics[width=4.2cm]{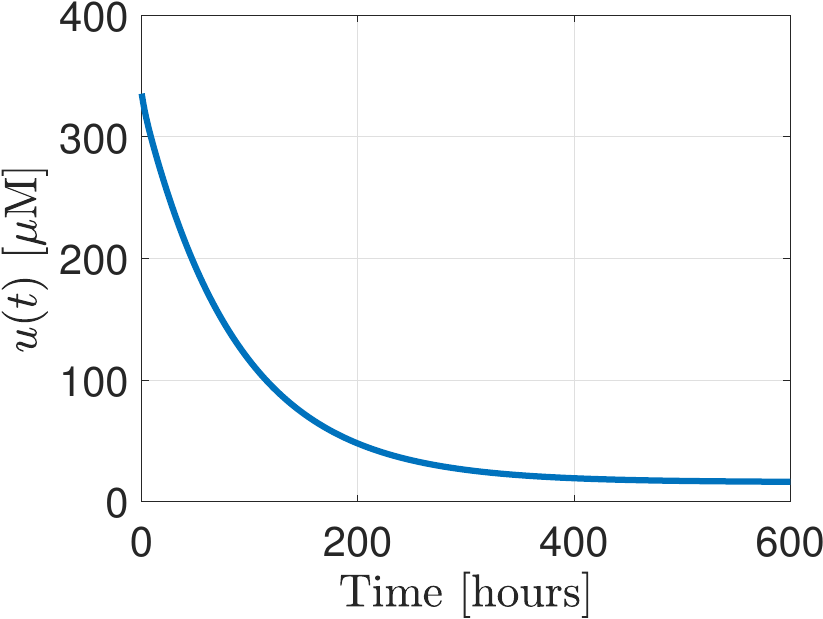}}
	\caption{Simulation results for the closed-loop system with the output-feedback control law (\ref{eq:increaseDosage}).}
	\label{fig:simulation}
	\begin{picture}(0,0)
		\put(-98,218){\includegraphics[width=2.75cm]{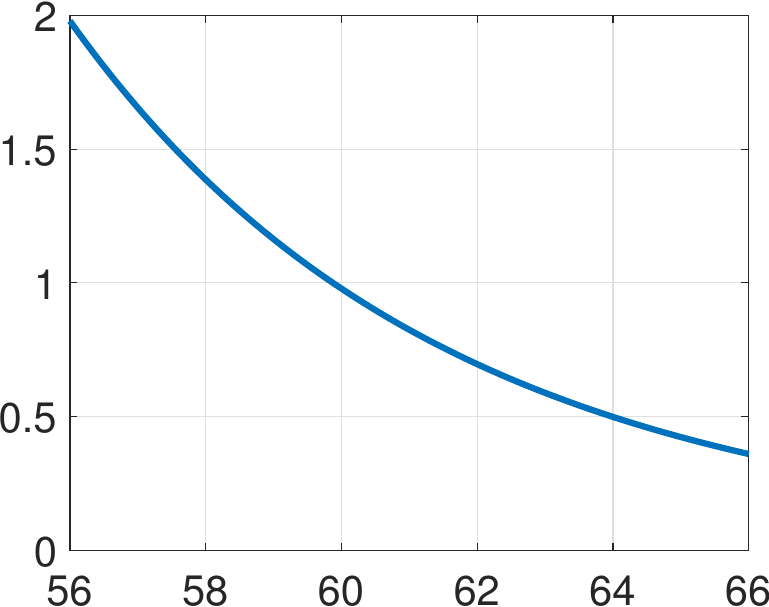}}
	\end{picture}
	\begin{picture}(0,0)
		\put(55,223){\includegraphics[width=2.25cm]{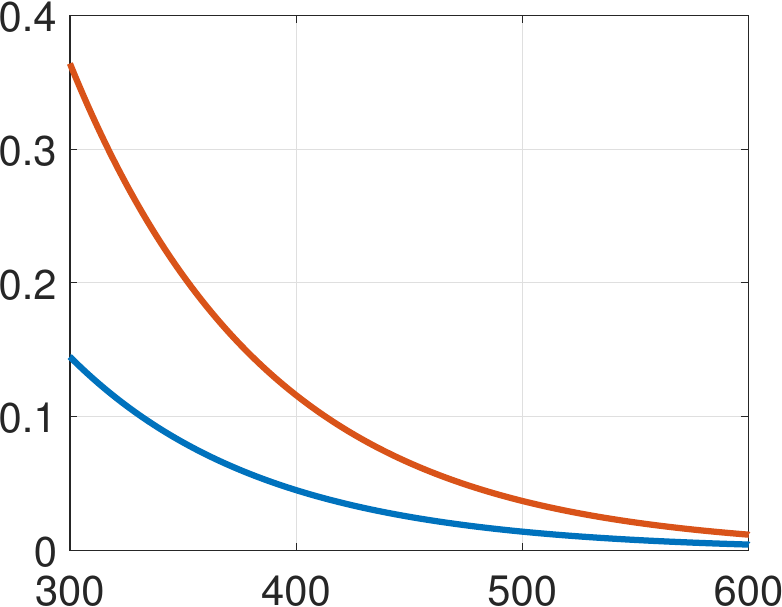}}
	\end{picture}	
		\begin{picture}(0,0)
		\put(-102,98){\includegraphics[width=2.75cm]{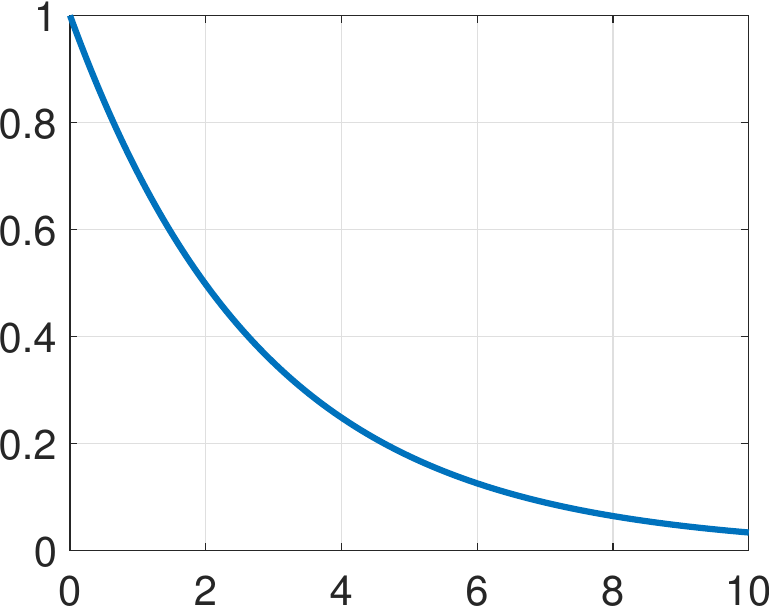}}
	\end{picture}	
		\begin{picture}(0,0)
		\put(45,108){\includegraphics[width=2.4cm]{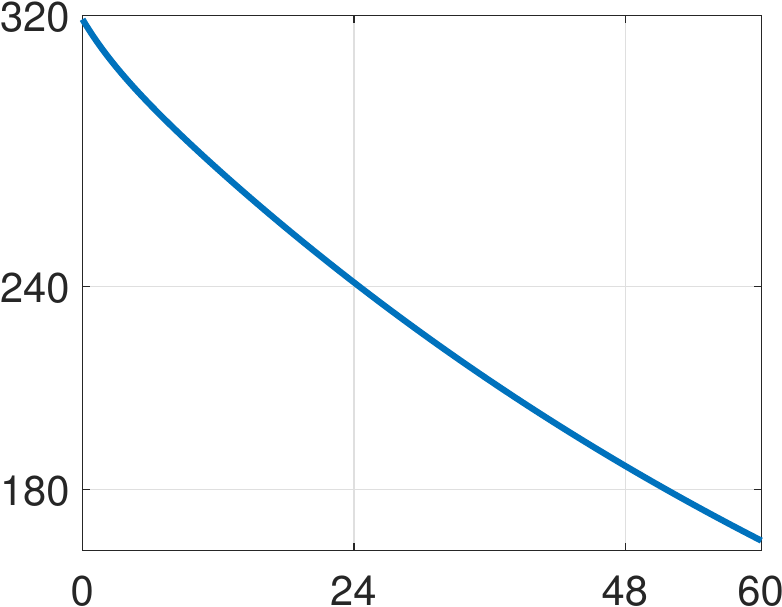}}
	\end{picture}		
\end{figure}

In Figure~\ref{fig:simulation}, the simulation results are presented with the proposed output-feedback control law (\ref{eq:increaseDosage}). Keeping the instantaneous action, Figure~\ref{fig:giardiaSimulation_y} shows that, after 60 hours from the beginning of the treatment, theoretically, the concentration of the protozoan is already lower than unity. In practice, for this ideal case, the giardia population will be null after 60h. Figure~\ref{fig:giardiaSimulation_yY0} shows that, in this situation, 10 hours after starting the treatment, the protozoan population is already less than 10\% of its initial value. Figure~\ref{fig:giardiaSimulation_x2hatX2} presents the norm observer behavior versus the unmeasured state variable. Figure~\ref{fig:giardiaSimulation_u} shows the decay of drug dosage over time generated by the control law (\ref{eq:increaseDosage}).

\section{Experiments}

\subsection{Axenic cultivation of {\it Giardia duodenalis} trophozoites and exposure to metronidazole} 

\begin{figure}[!ht]
	\centering
	\includegraphics[width=7.4cm]{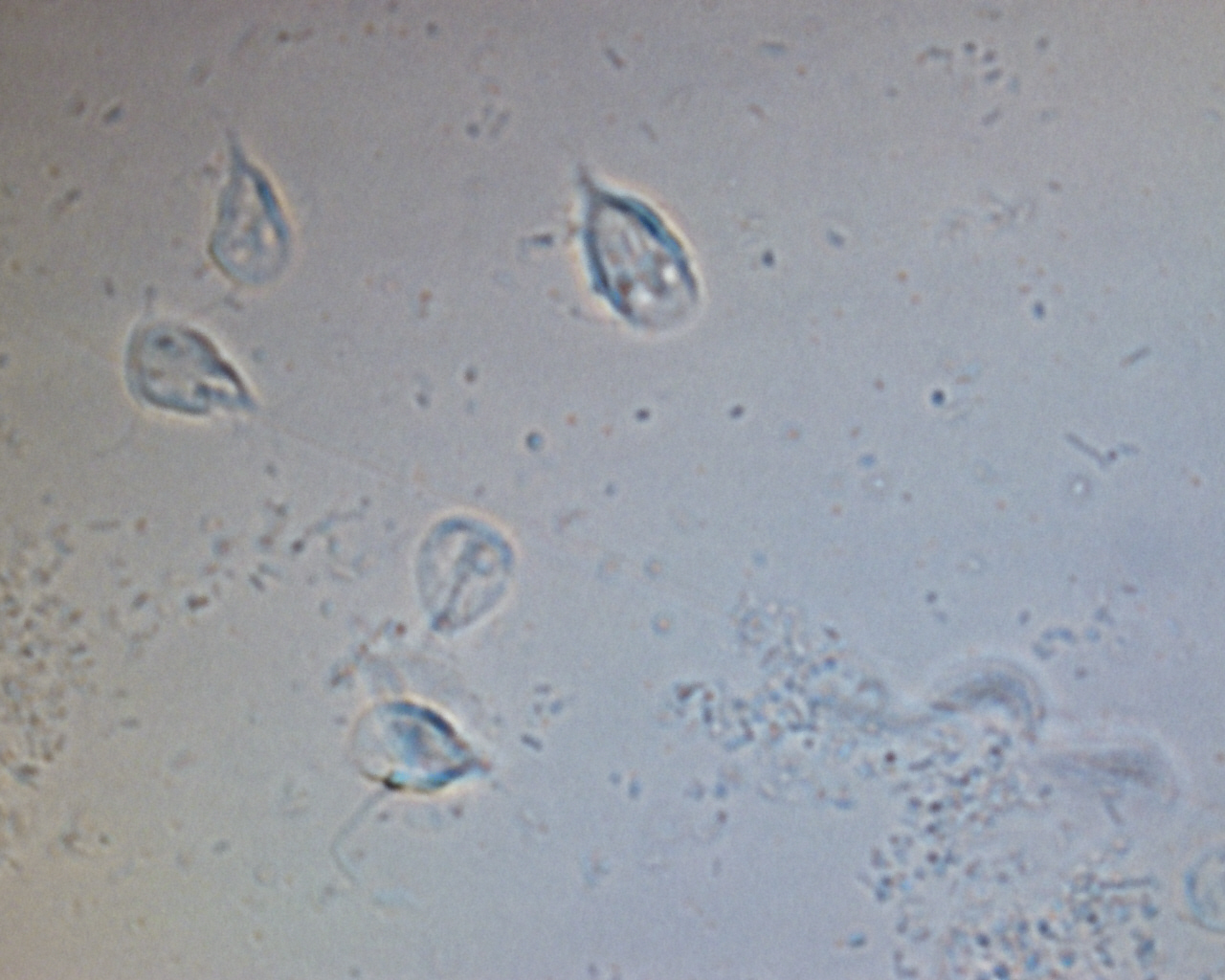}
	\caption{{\it Giardia lamblia} trophozoite by microscopy. \label{fig:giardiaExp}}
\end{figure}

In order to carry out the experimental study, the cultivation of the {\it Giardia} protozoan was first established. Trofozoites of the {\it G. duodenalis} species of clone W6, strain WB [ATCC50803], were grown in their own medium (TYI-S-33) \citep{K:1983}, at pH 7, supplemented with 10\% inactivated bovine serum, in tubes of centrifuge with a capacity of 15 mL, see Figure~\ref{fig:giardiaExp}. During the stationary phase of growth of the parasite (107 cells/mL), the trophozoites passages were performed to maintain the cells. The tubes containing the parasite were centrifuged at 840 g for 5 minutes at 25$^{o}$C, the supernatant was discarded and 10 ml of culture medium was added. To inhibit the adsorption of the parasites, the tubes were incubated on ice for 20 minutes and shaken with the aid of a vortex. A total of 500 $\mu$L of original culture was added in three centrifuge tubes containing new culture medium with 10\% fetal bovine serum with a final volume of 10 mL. After raising, the tubes were incubated in an oven at 37$^{o}$C.

After the establishment of the cultivation, four independent experiments were carried out simultaneously. A concentration of 106 cells per ml was exposed to different concentrations of the drug metronidazole (Sigma Chemical, USA) and followed up to 60 h. Every 24h period, cells were counted in a Neubauer improved bright-line and the number of cells in each experiment was estimated. At the time of counting, respecting the ideal conditions for the cultivation of trophozoites (culture medium and fetal bovine serum), a new exposure to the drug was carried out. The drug concentrations used in each experiment at the respective analysis times are shown in Table~1.

\begin{table}[!ht]
\begin{center}
\caption{Metronidazole concentration used to exposition of Giardia duodenalis trophozoite in each time.}
\begin{tabular}{|c|c|c|c|}
\hline
              &     0h     &   24h      &    48h     \\
\hline 
Experiment 1  & 320 $\mu$M & 240 $\mu$M & 180 $\mu$M \\
\hline
Experiment 2  & 320 $\mu$M & 240 $\mu$M & 160 $\mu$M \\
\hline
Experiment 3  & 320 $\mu$M & 160 $\mu$M &  80 $\mu$M \\
\hline
Experiment 4  & 320 $\mu$M &  80 $\mu$M &  40 $\mu$M \\
\hline
\end{tabular}
\end{center}
\end{table}

\begin{figure}[H]
	\centering
	\subfigure[{\it Giardia lamblia} Population. ]{\includegraphics[width=4.2cm]{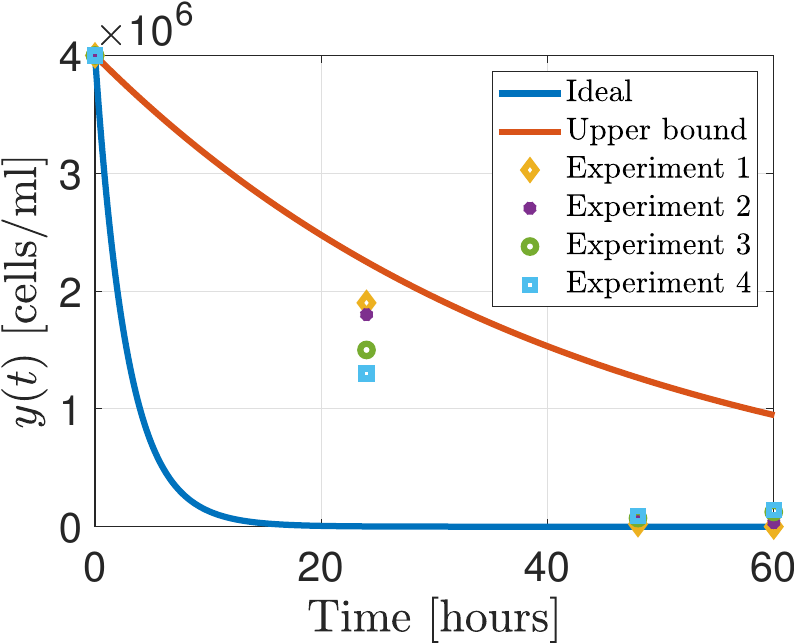}}
	~~
\subfigure[Drug dosage in micro molar ($1\mu g/ml\approx 5,8425\mu M$).]{\includegraphics[width=4.2cm]{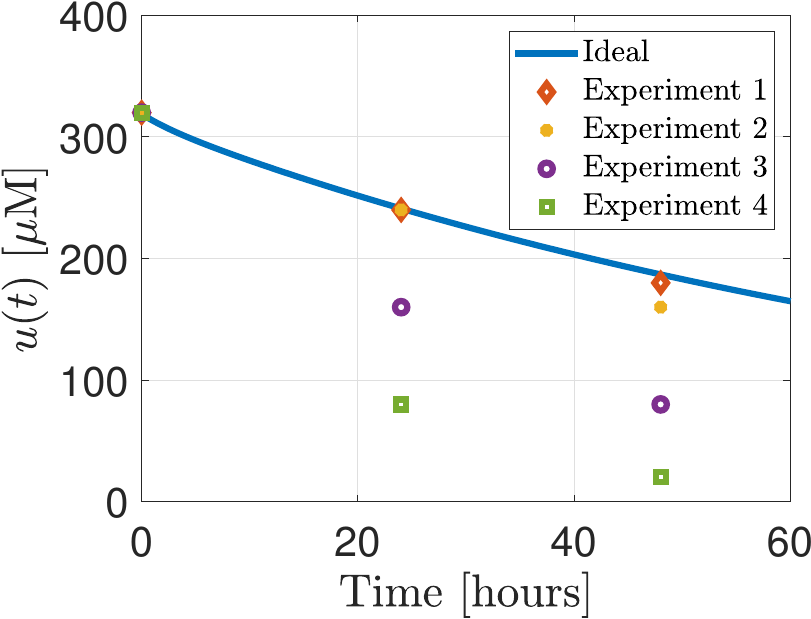}}
	\caption{Experiment results based on the output-feedback control law (\ref{eq:increaseDosage}).}
	\label{fig:experiment}
\end{figure}

The experimental results are shown in Figure~\ref{fig:experiment}. 

\subsection{Discussion}

As discussed before, Giardiasis is an infectious-parasitic disease that has an incidence of approximately 280 million cases each year \citep{LL:2002,L:2010}. The treatment of choice for the treatment of giardiasis consists of the use of 5-nitroimidazole drugs, mainly metronidazole \citep{GH:2001}. Treatment of the disease with metronidazole is generally effective, with cure rates approaching 90\% \citep{YE:2011}. However, the increase in cases of parasitic persistence after treatment has been reported in developed and developing countries \citep{NLAC:2015,FOFTVBC:2020}. The refractory cases observed can be the result of reinfection, parasitic resistance or even an inadequate dose of the medication.

In this study, extremely high exposure doses were used when compared to the doses used in clinical treatment. In the experiments, a higher initial dose (320 $\mu$M) was used and within 24 hours the dose was reduced (although still considered high). The methodology rationale predicted a greater impact in reducing the parasite number in the culture at the initial moment and providing smaller doses later would be able to eliminate exponentially the presence of trophozoites. 

In order to simulate what occurs naturally in the human intestine, throughout the experiment the trophozoites were not deprived of nutrients. Using this strategy, it was possible to exclude the possibility of parasite deaths due to the lack of ideal conditions and, thus, we could analyze the exclusive effect of the drug.

\section{Conclusion}

In this paper, unlike what was done in \citep{LLZ:2013} where all parameters were considered to be known and the full-state vector of the plant was available, the population control problem of {\it G. lamblia} was addressed considering that the model parameters are uncertain and only the output variable is available for feedback. By employing a norm observer, it was possible to find an upper bound for the absolute value of the unmeasured state variable that represents the mutation of the protozoan. Thus, by using known parametric bounds, it was possible to design an output-feedback control law. Simulation results were presented to illustrate the behavior of the system in closed loop and experimental tests were be carried out in order to experimentally validate the proposed control law.

Since the stabilization of the {\it Giardia lamblia} occours exponentially, even though a significant reduction in trophozoites was observed when exposed to high concentrations of metronidazole, as already reported in other studies \citep{OFC:2020}, a high number of cells was still ultimately observed after $60$ hours. It should be noted that although the cell counting methodology in the Neubauer has its limitations, and the scale of magnitude observed was still very high.

The findings of this study point out to an alarming data, since even very high doses of metronidazole can be inefficient to completely eliminate the parasite, in the 60h interval. As this protozoan reproduces by binary division, if not eliminated, it can proliferate again and colonize the region of the duodenum. This fact is worrying, since it can lead to prolonged contact of the parasite with the small intestine, leading to a picture of chronic giardiasis.

Extremum seeking control approaches for population control of {\it Giardia lamblia} using sampled-data \citep{Emilia_TAC} or sliding mode-based implementations
\citep{AOH:13,c5}, considering delays \citep{c22} and delay mismatch \citep{Damir_EJC}, wider classes of PDE models \citep{RAD_PDE_ALCOS,OFKK2020,PDE_cascades_SCL,Oliveira_Krstic_2022}, with dynamic maps \citep{Damir_LCSS} as well as game theory \citep{Basar_heat,Basar_heterogeneous,ORKBEK:19} are also possible topics for future investigation. For instance, an automatic control parameter tuning could be envisaged from the methodology introduced for neuromuscular electrical stimulation in \citep{POPF2020}.

\begin{ack}
This work was financed in part by the Coordena{\c c}{\~a}o de Aperfeiçoamento de Pessoal de N{\'i}vel Superior – Brasil (CAPES) – Finance Code 001. The authors also acknowledge the Brazilian Funding Agencies Conselho Nacional de Desenvolvimento Cient{\'i}fico e Tecnol{\'o}gico (CNPq) and Funda{\c c}{\~a}o de Amparo {\`a} Pesquisa do Estado do Rio de Janeiro (FAPERJ). This work was supported by CNPq Universal Program (Grant 435015/2018-4), FAPERJ (Grant E-26/202.078/2020), and Instituto Oswaldo Cruz/FIOCRUZ—Brazilian Ministério da Saúde (internal funds). The second author was supported by a fellowship CNPq (PDJ), INOVA FIOCRUZ Program and FAPERJ Nota 10.
\end{ack}


\begin{thebibliography}{xx}  

\bibitem[Aminde et~al.(2013)Aminde, Oliveira, and Hsu]{AOH:13}
N.O. Aminde, T.R. Oliveira, and L. Hsu.
\newblock Global output-feedback extremum seeking control via monitoring functions.
\newblock In \emph{52nd IEEE Conference on Decision and Control (CDC)}, pages 1031--1036, 2013.

\bibitem[Baldursson and Karanis(2011)]{BK:2011}
S. Baldursson and P. Karanis.
\newblock Waterborne transmission of protozoan parasites: Review of worldwide outbreaks -- An update 2004-2010.
\newblock \emph{Water Research}, 45:6603--6614, 2011.

\bibitem[Berkman et~al.(2002)Berkman, Lescano, Gilman, Lopez, and Black]{BLGLB:2002}
D.S. Berkman, A.G. Lescano, R.H. Gilman, S.L. Lopez, and M.M. Black.
\newblock Effects of stunting, diarrhoeal disease, and parasitic infection during infancy on cognition in late childhood: a follow-up study.
\newblock \emph{The Lancet}, 359:564--571, 2002.

\bibitem[Boreham et~al.(1988)Boreham, Phillips, and Shepherd]{BPS:1988}
P.F.L. Boreham, R.E. Phillips, and R.W. Shepherd.
\newblock Altered uptake of metronidazole in vitro by stocks of Giardia intestinalis with different drug sensitivities.
\newblock \emph{Transactions of the Royal Society of Tropical Medicine and Hygiene}, 82:104--106, 1988.

\bibitem[Carranza and Lujan(2010)]{CL:2010}
P.G. Carranza and H.D. Lujan.
\newblock New insights regarding the biology of Giardia lamblia.
\newblock \emph{Microbes and Infection}, 12:71--80, 2010.

\bibitem[Cunha et~al.(2008)Cunha, Costa, and Hsu]{CCH:2008}
J.P.V.S. Cunha, R.R. Costa, and L. Hsu.
\newblock Design of first-order approximation filters for sliding-mode control of uncertain systems.
\newblock \emph{IEEE Transactions on Industrial Electronics}, 55:4037--4046, 2008.

\bibitem[Daumerie and Kindhauser(2004)]{DK:2004}
D. Daumerie and M.K. Kindhauser.
\newblock Intensified Control of Neglected Diseases.
\newblock World Health Organization, 2004.

\bibitem[Fantinatti et~al.(2020)Fantinatti, Lopes-Oliveira, Cascais-Figueredo, Austriaco-Teixeira, Verissimo, Bello, and Da-Cruz]{FOFTVBC:2020}
M. Fantinatti, L.A.P. Lopes-Oliveira, T. Cascais-Figueredo, P. Austriaco-Teixeira, E. Verissimo, A.R. Bello, and A.M. Da-Cruz.
\newblock Recirculation of Giardia lamblia Assemblage A After Metronidazole Treatment in an Area With Assemblages A, B, and E Sympatric Circulation.
\newblock \emph{Frontiers in Microbiology}, 11:23--39, 2020.

\bibitem[Farthing(1996)]{F:1996}
M.J.G. Farthing.
\newblock Giardiasis.
\newblock \emph{Gastroenterology Clinics of North America}, 25:493--515, 1996.

\bibitem[Feng and Xiao(2011)]{FX:2011}
Y. Feng and L. Xiao.
\newblock Zoonotic potential and molecular epidemiology of Giardia species and giardiasis.
\newblock \emph{Clinical Microbiology Reviews}, 24:110--140, 2011.

\bibitem[Gardner and Hill(2001)]{GH:2001}
T.B. Gardner and D.R. Hill.
\newblock Treatment of Giardiasis.
\newblock \emph{Clinical Microbiology Reviews}, 14:114--128, 2001.

\bibitem[Guimarães and Sogayar(2002)]{GS:2002}
S. Guimarães and M.I.L. Sogayar.
\newblock Detection of anti-Giardia lamblia serum antibody among children of day care centers.
\newblock \emph{Revista de Saúde Pública}, 36:63--68, 2002.

\bibitem[Keister(1983)]{K:1983}
D.B. Keister.
\newblock Axenic culture of Giardia lamblia in TYI-S-33 medium supplemented with bile.
\newblock \emph{Transactions of The Royal Society of Tropical Medicine and Hygiene}, 77:487--488, 1983.

\bibitem[Khalil(2002)]{K:2002}
H.K. Khalil.
\newblock \emph{Nonlinear Systems}.
\newblock Prentice Hall, Upper Saddle River, New Jersey, 2002.

\bibitem[Lalle(2010)]{L:2010}
M. Lalle.
\newblock Giardiasis in the post genomic era: treatment, drug resistance and novel therapeutic perspectives.
\newblock \emph{Infectious Disorders - Drug Targets}, 10:283--294, 2010.

\bibitem[Lane and Lloyd(2002)]{LL:2002}
S. Lane and D. Lloyd.
\newblock Current Trends in Research into the Waterborne Parasite Giardia.
\newblock \emph{Critical Reviews in Microbiology}, 28:123--147, 2002.

\bibitem[Li et~al.(2013)Li, Lenaghan, and Zhang]{LLZ:2013}
X. Li, S.C. Lenaghan, and M. Zhang.
\newblock Evolutionary game based control for biological systems with applications in drug delivery.
\newblock \emph{Journal of Theoretical Biology}, 326:58--69, 2013.

\bibitem[Lima Junior et~al.(2013)Lima Junior, Kaiser, and Catisti]{JKC:2013}
O.A. Lima Junior, J. Kaiser, and R. Catisti.
\newblock High occurrence of giardiasis in children living on a ``landless farm workers'' settlement in Araras, São Paulo, Brazil.
\newblock \emph{Revista do Instituto de Medicina Tropical de Saúde Pública}, 55:185--188, 2013.

\bibitem[Lopes-Oliveira et~al.(2020)Lopes-Oliveira, Fantinatti, and Da-Cruz]{OFC:2020}
L.A.P. Lopes-Oliveira, M. Fantinatti, and A.M. Da-Cruz.
\newblock In vitro-induction of metronidazole-resistant Giardia duodenalis is not associated with nucleotide alterations in the genes involved in pro-drug activation.
\newblock \emph{Memórias do Instituto Oswaldo Cruz}, 115:1--4, 2020.

\bibitem[Nabarro et~al.(2015)Nabarro, Lever, Armstrong, and Chiodini]{NLAC:2015}
L.E.B. Nabarro, R.A. Lever, M. Armstrong, and P.L. Chiodini.
\newblock Increased incidence of nitroimidazole-refractory giardiasis at the Hospital for Tropical Diseases, London: 2008--2013.
\newblock \emph{Clinical Microbiology and Infection}, 21:791--796, 2015.

\bibitem[Oliveira et~al.(2019)Oliveira, Feiling, and Krstic]{RAD_PDE_ALCOS}
T.R. Oliveira, J. Feiling, and M. Krstic.
\newblock Extremum seeking for maximizing higher derivatives of unknown maps in cascade with reaction-advection-diffusion PDEs.
\newblock \emph{IFAC-PapersOnLine}, 52:210--215, 2019.

\bibitem[Oliveira et~al.(2020)Oliveira, Feiling, Koga, and Krstic]{OFKK2020}
T.R. Oliveira, J. Feiling, S. Koga, and M. Krstic.
\newblock Extremum seeking for unknown scalar maps in cascade with a class of parabolic partial differential equations.
\newblock \emph{International Journal of Adaptive Control and Signal Processing}, 35(7):1162--1187, 2020.

\bibitem[Oliveira et~al.(2011)Oliveira, Hsu, and Peixoto]{c5}
T.R. Oliveira, L. Hsu, and A.J. Peixoto.
\newblock Output-feedback global tracking for unknown control direction plants with applications to extremum-seeking control.
\newblock \emph{Automatica}, 47:2029--2038, 2011.

\bibitem[Oliveira and Krstic(2015)]{c22}
T.R. Oliveira and M. Krstic.
\newblock Newton-based extremum seeking under actuator and sensor delays.
\newblock \emph{IFAC-PapersOnLine}, 48:304--309, 2015.

\bibitem[Oliveira and Krstic(2021)]{PDE_cascades_SCL}
T.R. Oliveira and M. Krstic.
\newblock Extremum seeking boundary control for PDE-PDE cascades.
\newblock \emph{Systems \& Control Letters}, 155:105004/1--105004/15, 2021.

\bibitem[Oliveira and Krstic(2022)]{Oliveira_Krstic_2022}
T.R. Oliveira and M. Krstic.
\newblock \emph{Extremum Seeking through Delays and PDEs}.
\newblock SIAM, Philadelphia, USA, 2022.

\bibitem[Oliveira et~al.(2021a)Oliveira, Rodrigues, Krstic, and Basar]{Basar_heat}
T.R. Oliveira, V.H.P. Rodrigues, M. Krstic, and T. Basar.
\newblock Nash equilibrium seeking with players acting through heat PDE dynamics.
\newblock In \emph{American Control Conference (ACC)}, pages 684--689, 2021a.

\bibitem[Oliveira et~al.(2021b)Oliveira, Rodrigues, Krstic, and Basar]{Basar_heterogeneous}
T.R. Oliveira, V.H.P. Rodrigues, M. Krstic, and T. Basar.
\newblock Nash equilibrium seeking in heterogeneous noncooperative games with players acting through heat PDE dynamics and delays.
\newblock In \emph{60th IEEE Conference on Decision and Control (CDC)}, pages 1167--1173, 2021b.

\bibitem[Oliveira et~al.(2021c)Oliveira, Rodrigues, Krstic, and Basar]{ORKBEK:19}
T.R. Oliveira, V.H.P. Rodrigues, M. Krstic, and T. Basar.
\newblock Nash equilibrium seeking in quadratic noncooperative games under two delayed information-sharing schemes.
\newblock \emph{J Optim Theory Appl}, 191:700--735, 2021c.

\bibitem[Paz et~al.(2020)Paz, Oliveira, Pino, and Fontana]{POPF2020}
P. Paz, T.R. Oliveira, A.V. Pino, and A.P. Fontana.
\newblock Model-free neuromuscular electrical stimulation by stochastic extremum seeking.
\newblock \emph{IEEE Transactions on Control Systems Technology}, 28(1):238--253, 2020.

\bibitem[Rusiti et~al.(2019)Rusiti, Evangelisti, Oliveira, Gerdts, and Krstic]{Damir_LCSS}
D. Rusiti, G. Evangelisti, T.R. Oliveira, M. Gerdts, and M. Krstic.
\newblock Stochastic extremum seeking for dynamic maps with delays.
\newblock \emph{IEEE Control Systems Letters}, 3:61--66, 2019.

\bibitem[Rusiti et~al.(2021)Rusiti, Oliveira, Krstic, and Gerdts]{Damir_EJC}
D. Rusiti, T.R. Oliveira, M. Krstic, and M. Gerdts.
\newblock Robustness to delay mismatch in extremum seeking.
\newblock \emph{European Journal of Control}, 62:75--83, 2021.

\bibitem[Tejman-Yarden and Eckmann(2011)]{YE:2011}
N. Tejman-Yarden and L. Eckmann.
\newblock New approaches to the treatment of giardiasis.
\newblock \emph{Current Opinion in Infectious Diseases}, 24:451--456, 2011.

\bibitem[Thompson et~al.(1993)Thompson, Reynoldson, and Mendis]{TRM:2013}
R.C. Thompson, J.A. Reynoldson, and A.H. Mendis.
\newblock Giardia and giardiasis.
\newblock In \emph{Advances in Parasitology}, volume 32, pages 71--160, 1993.

\bibitem[Zaat et~al.(1997)Zaat, Mank, and Assendelft]{ZMA:1997}
J.O.M. Zaat, T.G. Mank, and W.J.J. Assendelft.
\newblock A systematic review on the treatment of giardiasis.
\newblock \emph{Tropical Medicine \& International Health}, 2:63--82, 1997.

\bibitem[Zhu et~al.(2023)Zhu, Fridman, and Oliveira]{Emilia_TAC}
Y. Zhu, E. Fridman, and T.R. Oliveira.
\newblock Sampled-data extremum seeking with constant delay: a time-delay approach.
\newblock \emph{IEEE Transactions on Automatic Control}, 68:432--439, 2023.

\end{thebibliography}

\end{document}